\documentclass[12pt]{article}
\usepackage{amsxtra}
\usepackage{amssymb, amsmath}
\usepackage{amsthm}
\usepackage[latin1]{inputenc}
\usepackage {hyperref}

\usepackage{amssymb}
\usepackage{amsmath}
\usepackage[dvips]{graphicx}

\begin{document}
\begin{center}
\textbf{\LARGE{Bimagic Squares of Bimagic Squares \\and an Open Problem}}
\end{center}

\bigskip
\begin{center}
\textbf{\large{Inder Jeet Taneja}}\\
Departamento de Matem\'{a}tica\\
Universidade Federal de Santa Catarina\\
88.040-900 Florian\'{o}polis, SC, Brazil.\\
\textit{e-mail: taneja@mtm.ufsc.br\\
http://www.mtm.ufsc.br/$\sim $taneja}
\end{center}

\begin{abstract}
\textit{In this paper we have produced different kinds of bimagic squares based on bimagic squares of order 8}$\times $\textit{8, 16}$\times $\textit{16, 25}$\times $\textit{25, 49}$\times $\textit{49, etc. A different technique is applied to produce bimagic square of order 16}$\times $\textit{16, 25}$\times $\textit{25, 49}$\times $\textit{49, etc. The bimagic square of order 8}$\times $\textit{8 used is the already known in the literature. The work is neither based on any programming language nor on mathematical results. Just simple combinations are used to produce these bimagic squares. Moreover, in each case we have used consecutive numbers starting from 1.}
\end{abstract}

\section{Introduction}

In this work we produced different orders of bimagic squares containing bimagic or semi-bimagic squares. Most of the work is based on the bimagic squares of orders 8$\times $8, 16$\times $16 and 25$\times $25. We have also produced new bimagic squares of order 49$\times $49, 121$\times $121. The bimagic square of order 8$\times $8 is obtained long back by Pfeffermann 1891 \cite{pfe}.  The bimagic squares of order 16$\times $16 and 25$\times $25 are very much similar to one given in \cite{dee} as examples. The difference is that the bimagic square of order 16$\times $16 appearing in \cite{dee} don't have the property that each sub-block of order 4$\times $4 as a magic square, while in our case it happens. The bimagic square of order 25$\times $25 is very much similar to given in \cite{dee}  but we have applied a little different approach to produce it. The bimagic square of order 49$\times $49 is long back produced by G. Tarry 1895 \cite{boy}. Here also we applied a little different approach. We have produced the bimagic square of order 121$\times $121 without knowing that is it done before. Based on the approach adopted in this work, we can always produce bimagic squares using squares of prime number such as \textbf{$13^{2} \times 13^{2} $}, \textbf{$17^{2} \times 17^{2} $}, etc.

 \bigskip
No programming language is used, just simple combinations are sufficient to produce whole the work. In each case, we have used consecutive numbers starting from 1. Some of these files are available at the authors' web-site given above.

\bigskip
During construction we observe that in case of orders $k^{4} \times k^{4} ,\; k=4,\; 5$ and $7$, in the previous subgroup $k^{3} \times k^{3} ,\; k=4,\; 5$ and  $7$ the bimagic sums has the value in each case. For more details see the table given at the end as an open problem to prove it mathematically.

\bigskip
Before we proceed, here below are some basic definitions:

\begin{itemize}
\item [(i)] A \textbf{magic square} is a collection of numbers put as a square matrix, where the sum of element of each row, sum of element of each column and sum of each element of two principal diagonals have the same sum.  For simplicity, let us write it as \textbf{S1}.
\item  [(ii)] \textbf{Bimagic square} is a magic square where the sum of square of each element of rows, columns and two principal diagonals are the same. For simplicity, let us write it as \textbf{S2}.
\item  [(iii)] \textbf{Upside down}, i.e., if we rotate it to $180^{0}$ degree it remains the same.
\item  [(iv)] \textbf{Mirror looking}, i.e., if we put it in front of mirror or see from the other side of the glass, or see on the other side of the paper, it always remains the magic square.
\item [(v)] \textbf{Universal magic squares}, i.e., magic squares having the property of upside down and mirror looking are considered \textit{universal magic squares}.
\end{itemize}

A good collection of multimagic squares can be seen in \cite{boy}. New upside down and universal magic squares can be seen in Taneja \cite{tan1}-\cite{tan7}.

\section{Details}

Whole the work we have divided in small parts. In the end we have given magic and bimagic sums in each case.

\subsection{First Part}

In this part we have presented bimagic squares of the following orders:

\begin{center}
16$\times $16, 32$\times $32, 56$\times $56, 64$\times $64, 72$\times $72, 88$\times $88, 96$\times $96, \\104$\times $104, 112$\times $112, 128$\times $128, 144$\times $144,  176$\times $176, 208$\times $208,\\ 224$\times $224, 256$\times $256, 512$\times $512, 1024$\times $1024, 2048$\times $2048 and 4096$\times $4096.
\end{center}

Let us divide the above bimagic squares in three small groups.

\subsubsection{First Small Group}

In this subsection we have given the following bimagic squares
\begin{center}
16$\times $16, 64$\times $64, 256$\times $256, 1024$\times $1024 and 4096$\times $4096.
\end{center}

In this group we have the special property that each block of orders 16$\times $16, 64$\times $64, 256$\times $256 and 1024$\times $1024  are also bimagic squares. Also, each block of order 4$\times $4 is a magic square. In case of magic squares of order 256$\times $256 each block of order 64$\times $64 has the same bimagic sum S2. In case of magic squares of order 1024$\times $1024 and 4096$\times $4096 each block of order 256$\times $256 produces a same bimagic sum S2 for 64$\times $64.
Using the same procedure we can also calculate the bimagic square of order 128$\times $128. Its values are given in the last section.

\subsubsection{Second Small Group}

In this subsection we have given the following bimagic squares

\begin{center}
32$\times $32, 64$\times $64, 96$\times $96, 128$\times $128, 512$\times $512 and 2048$\times $2048.
\end{center}

This group has the special property that each block of orders 8$\times $8, 64$\times $64, 96$\times $96, 128$\times $128, 512$\times $512 and 2048$\times $2048 is either a bimagic or semi-bimagic square but the final group is always a bimagic square. Each block of order 8$\times $8 is always a magic square having the same magic sum S1 in whole the order.
The bimagic square of order 160$\times $160 can also be calculated with the same procedure, but we have calculated it in the next part as multiple of 16. If we calculate 128$\times $128 according to first small group i.e., using bimagic square of 16$\times $16, then the next orders 512$\times $512 and 2048$\times $2048 can also be calculated as combinations of bimagic squares of order 16$\times $16. In this case all the bimagic squares of this group goes to first small group except 32$\times $32 and 96$\times $96.

\subsubsection{Third Small Group}

In this subsection we have given the following bimagic squares
\begin{center}
 56$\times $56, 72$\times $72, 88$\times $88, 104$\times $104, 112$\times $112,  \\ 144$\times $144,  176$\times $176, 208$\times $208  and 224$\times $224.
\end{center}

In case of 56$\times $56, 72$\times $72, 88$\times $88 and  104$\times $104, each block of order  8$\times $8 is either a bimagic or semi-bimagic square but the final order is always a bimagic square. Each block of order 8$\times $8 is always a magic square. While, in case of 112$\times $112, 144$\times $144,   176$\times $176, 208$\times $208 and  224$\times $224 each block of order 16$\times $16 is a bimagic square with the property that each block of order 4$\times $4 is a magic magic square.

\subsection{Second Part}

In this part we have presented bimagic squares of the following orders:

\begin{center}
40$\times $40, 80$\times $80, 120$\times $120, 160$\times $160, 200$\times $200, 240$\times $240,  \\ 400$\times $400, 600$\times $600, 800$\times $800, 960$\times $960, 1000$\times $1000, 1200$\times $1200, \\ 1600$\times $1600,  2000$\times $2000, 2400$\times $2400, 3000$\times $3000, 3200$\times $3200 and 4000$\times $4000.
\end{center}

Let us divide the above bimagic squares in two small groups

\subsubsection{First Small Group}

In this subsection we have given the following bimagic squares

\begin{center}
80$\times $80, 160$\times $160,  240$\times $240, 400$\times $400, 800$\times $800, 960$\times $960,  \\
1200$\times $1200, 1600$\times $1600, 2000$\times $2000, 3200$\times $3200  and 4000$\times $4000.
\end{center}

This group we have the special property that each block of orders 16$\times $16 is a bimagic square with each block of order 4$\times $4 as a magic square. In higher cases such as 400$\times $400, each block of order 80$\times $80 is also a bimagic square, etc.

\subsubsection{Second Small Group}

\begin{center}
40$\times $40, 120$\times $120,  200$\times $200, 600$\times $600, \\ 1000$\times $1000,  2400$\times $2400 and 3000$\times $3000,
\end{center}

This group has the special property that each block of orders 8$\times $8 is either a bimagic or semi-bimagic but the final group is always a bimagic square. Each block of order 8$\times $8 is always a magic square. In the higher case such as 200$\times $200, each block of order 40$\times $40 is a bimagic or semi-bimagic, etc.

\subsection{Third Part}

In this section we have presented bimagic squares of following orders:

\begin{center}
25$\times $25, 125$\times $125 and 625$\times $625
\end{center}

Here each block of order 5$\times $5 is a magic square with the same sum S1. Each block of orders 25$\times $25  is a bimagic square. In case of  625$\times $625, each block of order 125$\times $125 is also a bimagic square with the same sum S2.

\subsection{ Forth Part}

In this section we have presented bimagic squares of the following orders:

\begin{center}
49$\times $49, 343$\times $343 and 2401$\times $2401.
\end{center}

Here each block of order 7$\times $7 is a magic square with the same sum S1. Each block of order 49$\times $49 is a bimagic square. Moreover, in case of  2401$\times $2401 each block of order 343$\times $343 has the same bimagic sum S2.

\subsection{Fifth Part}

In this section we have presented bimagic squares of the following orders:

\begin{center}
121$\times $121 and 1331$\times $1331.
\end{center}

Here each block of order 11$\times $11 is a magic square with the same sum S1. Each block of order 121$\times $121 is a bimagic square. The next  order 14641$\times $14641 multiple of 11 with 1331  is not calculated because of higher values.

\section{Semi-Bimagic Squares}

 Following the same approach applied to obtain bimagic squares in the above sections, we can still have semi-bimagic squares of orders 24$\times $24 and 48$\times $48 with the property that in case of 24$\times $24, each block of order 8$\times $8 is a magic square with three blocks of order 8$\times $8 as bimagic and six blocks as semi-bimagic. In case of 48$\times $48, each block of order 4$\times $4 is a magic square. Here we have 9 blocks of bimagic squares of order 16$\times $16. For the numerical values of these semi-bimagic squares see the last section.

 \section{Bimagic Squares of Magic Squares}

Tarry-Cazalas \cite{gec} in 1934 (see \cite{boy}) gave a bimagic square of order 9$\times $9 . Following the idea of Tarry-Cazalas and the approach adopted above, we obtained here bimagic squares of orders:

\begin{center}
 81$\times $81 and 729$\times $729 .
\end{center}

 These two bimagic squares have the property that each block of order 9$\times $9 is just a magic square, rathar than bimagic as in the other cases studied above. Also the sum of all 9 member in each block of order 3$\times $3 has the same sum as of $S1_{9\times 9}$.  If we follow the idea of Pfeffermann  (see \cite{boy}) of bimagic square of order 9$\times $9 and use our approach, we are unable to get bimagic squares of orders 81$\times $81 and 729$\times $729. In whole the work, this is the only case, where we don't have subgroups of bimagic squares.

 We can still have bimagic squares of orders 81$\times $81 and 729$\times $729  considering just magic square of order  9$\times $9. In this situation, we get bimagic squares with the property that each small groups of orders 3$\times $3 and  9$\times $9 are semi-magic squares finally giving bimagic squares of orders 81$\times $81 and 729$\times $729 . For details see the files given in authors' site \cite{tan9}.

\section{Open Problem}

The bimagic square of order 16$\times $16 already known in the literature don't have the property that each block of order 4$\times $4 as a magic square. This we have done in this work. We have produced bimagic squares of order 25$\times $25, 49$\times $49, 121$\times $121, etc. in a different approach than the one already known in the literature. Interesting the approach adopted lead us the following property:

\bigskip
\noindent
\begin{tabular}{|p{0.09in}|p{1.68in}|p{1.85in}|p{1.05in}|p{0.95in}|} \hline
 & Bimagic Squares & Sub group-1 & Sub group-2 & Sub group-3 \\ \hline
1. & \textbf{Order 256$\times $256\newline }S1:=8388736\newline S2:=366512264576 &\textbf{Order 64$\times $64}\newline S1:= 2097184\newline \textbf{S2:=91628066144} & \textbf{Order16$\times $16}\newline S1:=524296\newline S2 is different in each case. & \textbf{Order 4$\times $4}\newline S1:=131074 \\ \hline
2. & \textbf{Order 625$\times $625\newline }S1:=122070625\newline S2:=31789265950625 & \textbf{Order 125$\times $125}\newline S1:=24414125\newline \textbf{S2:=6357853190125} & \textbf{Order 25$\times $25}\newline S1:=4882825\newline S2 is different in each case. & \textbf{Order 5$\times $5} \newline S1:=976565 \\ \hline
3. & \textbf{Order 2401$\times $2401\newline }S1:=6290644801\newline S2:=26597429019848000 & \textbf{Order 343$\times $343}\newline S1:=988663543\newline \textbf{S2:=3799632717121140} & \textbf{Order 49$\times $49}\newline S1:=141237649\newline S2 is different in each case. & \textbf{Order 7$\times $7}\newline S1:=20176807 \\ \hline
\end{tabular}

\bigskip
We observe from the above table that in the first sub-group in each case the S2 is same, i.e., $k^{4} \to k^{3} \to k^{2} \to k^{1} ,  \; k=4,\; 5$ and 7. Now the question is to prove it mathematically? Moreover, if we go to higher order, for example in case of 1024$\times $1024, the bimagic squares of order 64$\times $64 has the same sum S2 in each group of 256$\times $256, while the bimagic sums S2 for 256$\times $256 give different values.

\section{Numerical Values}

Here below are the numerical values in each case. Some files of these bimagic squares can be downloaded at authors' web-site \cite{tan9} :

\bigskip
\noindent
\textbf{$\bullet$ Bimagic square of order 16$\times $16}

\bigskip
$S1_{4\times 4} :=\frac{1}{4} S1_{16\times 16} =514$

$S1_{16\times 16} :=2056$

$S2_{16\times 16} :=351576$

\bigskip
\noindent
\textbf{$\bullet$ Semi-Bimagic square of order 24$\times $24}

\bigskip
$S1_{8\times 8 :=\frac{1}{3} S1_{24\times 24}} =2308$

$S1_{24\times 24} :=6924$

$S2_{24\times 24} :=2661124$ - (rows and columns)

$S2_{24\times 24} :=2654292$ - (diagonal - 1)

$S2_{24\times 24} :=2714116$ - (diagonal - 2)

\bigskip
\noindent
\textbf{$\bullet$ Bimagic square of order 25$\times $25}

\bigskip
$S1_{5\times 5} :=\frac{1}{5} S1_{25\times 25} =1565$

$S1_{25\times 25} :=7825$

$S2_{25\times 25} :=3263025$

\bigskip
\noindent
\textbf{$\bullet$ Bimagic square of order 32$\times $32}

\bigskip
$S1_{8\times 8} :=\frac{1}{4} S1_{32\times 32} =4100$

$S1_{32\times 32} :=16400$

$S2_{32\times 32} :=11201200$

\bigskip
\noindent
\textbf{$\bullet$ Bimagic square of order 40$\times $40}

\bigskip
$S1_{8\times 8} :=\frac{1}{5} S1_{40\times 40} =6404$

$S1_{40\times 40} :=32020$

$S2_{40\times 40} :=32165340$

\bigskip
\noindent
\textbf{$\bullet$ Semi-Bimagic square of order 48$\times $48}

\bigskip
$S1_{4\times 4} :=\frac{1}{4} S1_{16\times 16} =\frac{1}{12} S1_{48\times 48} =4610$

$S1_{48\times 48} :=55320$

$S2_{48\times 48} :=84989960$ - (rows and columns)

$S2_{48\times 48} :=84990120$ - (diagonal - 1)

$S2_{48\times 48} :=85358600$ - (diagonal - 2)

\bigskip
\noindent
\textbf{$\bullet$ Bimagic square of order 49$\times $49}

\bigskip
$S1_{7\times 7} :=\frac{1}{7} S1_{49\times 49} =8407$

$S1_{49\times 49} :=58849$

$S2_{49\times 49} :=94217249$

\bigskip
\noindent
\textbf{$\bullet$ Bimagic square of order 56$\times $56}

\bigskip
$S1_{8\times 8} :=\frac{1}{7} S1_{56\times 56} =12548$

$S1_{56\times 56} :=87836$

$S2_{56\times 56} :=183665076$

\bigskip
\noindent
\textbf{$\bullet$ Bimagic square of order 64$\times $64}

\bigskip
$S1_{8\times 8} :=\frac{1}{8} S1_{64\times 64} =16338$

or

$S1_{4\times 4} :=\frac{1}{4} S1_{16\times 16} =\frac{1}{16} S1_{64\times 64} =8194$

$S1_{64\times 64} :=131104$

$S2_{64\times 64} :=358045024$

\bigskip
\noindent
\textbf{$\bullet$ Bimagic square of order 72$\times $72}

\bigskip
$S1_{8\times 8} :=\frac{1}{9} S1_{72\times 72} =20740$

$S1_{72\times 72} :=186660$

$S2_{72\times 72} :=645159180$

\bigskip
\noindent
\textbf{$\bullet$ Bimagic square of order 80$\times $80}

\bigskip
$S1_{4\times 4} :=\frac{1}{4} S1_{16\times 16} =\frac{1}{20} S1_{80\times 80} =12802$

$S1_{80\times 80} :=256040$

$S2_{80\times 80} :=1092522680$

\bigskip
\noindent
\textbf{$\bullet$ Bimagic square of order 81$\times $81}

\bigskip
$S1_{9\times 9} :=\frac{1}{9} S1_{81\times 81} =29529$

$S1_{81\times 81} :=442416$

$S2_{81\times 81} :=1162527201$

\bigskip
\noindent
\textbf{$\bullet$ Bimagic square of order 88$\times $88}

\bigskip
$S1_{8\times 8} =\frac{1}{11} S1_{88\times 88} =30980$

$S1_{88\times 88} :=340780$

$S2_{88\times 88} :=1759447140$

\bigskip
\noindent
\textbf{$\bullet$  Bimagic square of order 96$\times $96}

\bigskip
$S1_{8\times 8} :=\frac{1}{12} S1_{96\times 96} =36868$

$S1_{96\times 96} :=442416$

$S2_{120\times 120} :=2718351376$

\bigskip
\noindent
\textbf{$\bullet$ Bimagic square of order 104$\times $104}

\bigskip
$S1_{8\times 8} =\frac{1}{13} S1_{104\times 104} =43268$

$S1_{104\times 104} :=562484$

$S2_{104\times 104} :=4056072124$

\bigskip
\noindent
\textbf{$\bullet$  Bimagic square of order 112$\times $112}

\bigskip
$S1_{4\times 4} :=\frac{1}{4} S1_{16\times 16} =\frac{1}{28} S1_{112\times 112} =25090$

$S1_{112\times 112} :=702520$

$S2_{112\times 112} :=5875174760$

\bigskip
\noindent
\textbf{$\bullet$ Bimagic square of order 120$\times $120}

\bigskip
$S1_{8\times 8} :=\frac{1}{15} S1_{120\times 120} =57604$

$S1_{120\times 120} :=864060$

$S2_{120\times 120} :=8295264020$

\bigskip
\noindent
\textbf{$\bullet$ Bimagic square of order 121$\times $121}

\bigskip
$S1_{11\times 11} :=\frac{1}{11} S1_{121\times 121} =80531$

$S1_{121\times 121} :=885841$

$S2_{121\times 121} :=8646694001$

\bigskip
\noindent
\textbf{$\bullet$ Bimagic square of order 125$\times $125}

\bigskip
$S1_{5\times 5} :=\frac{1}{5} S1_{25\times 25} =\frac{1}{25} S1_{125\times 125} =39065$

$S1_{125\times 125} :=976625$

$S2_{125\times 125} :=10173502625$

\bigskip
\noindent
\textbf{$\bullet$ Bimagic square of order 128$\times $128}

\bigskip
$S1_{8\times 8} :=\frac{1}{4} S1_{32\times 32} =\frac{1}{16} S1_{128\times 128} =65540$

or

$S1_{4\times 4} :=\frac{1}{4} S1_{16\times 16} =\frac{1}{32} S1_{128\times 128} =32770$

$S1_{128\times 128} :=1048640$

$S2_{128\times 128} :=11454294720$

\bigskip
\noindent
\textbf{$\bullet$  Bimagic square of order 144$\times $144}

\bigskip
$S1_{4\times 4} :=\frac{1}{4} S1_{16\times 16} =\frac{1}{36} S1_{144\times 144} =41474$

$S1_{144\times 144} :=1493064$

$S2_{144\times 144} :=20640614424$

\bigskip
\noindent
\textbf{$\bullet$ Bimagic square of order 160$\times $160}

\bigskip
$S1_{4\times 4} :=\frac{1}{4} S1_{16\times 16} =\frac{1}{40} S1_{160\times 160} =51202$

$S1_{160\times 160} :=2048080$

$S2_{160\times 160} :=34954581360$

\bigskip
\noindent
\textbf{$\bullet$ Bimagic square of order 176$\times $176}

\bigskip
$S1_{4\times 4} :=\frac{1}{4} S1_{16\times 16} =\frac{1}{44} S1_{176\times 176} =61954$

$S1_{176\times 176} :=2725976$

$S2_{176\times 176} :=56294130376$

\bigskip
\noindent
\textbf{$\bullet$ Bimagic square of order 200$\times $200}

\bigskip
$S1_{8\times 8} :=\frac{1}{5} S1_{40\times 40} =\frac{1}{25} S1_{200\times 200} =160004$

$S1_{200\times 200} :=4000100$

$S2_{200\times 200} :=106670666700$

\bigskip
\noindent
\textbf{$\bullet$ Bimagic square of order 208$\times $208}

\bigskip
$S1_{4\times 4} :=\frac{1}{4} S1_{16\times 16} =\frac{1}{52} S1_{208\times 208} =86530$

$S1_{208\times 208} :=4499560$

$S2_{208\times 208} :=129780809080$

\bigskip
\noindent
\textbf{$\bullet$ Bimagic square of order 224$\times $224}

\bigskip
$S1_{4\times 4} :=\frac{1}{4} S1_{16\times 16} =\frac{1}{52} S1_{224\times 224} =100354$

$S1_{224\times 224} :=5619824$

$S2_{224\times 224} :=187988732624$

\bigskip
\noindent
\textbf{$\bullet$ Bimagic square of order 240$\times $240}

\bigskip
$S1_{4\times 4} :=\frac{1}{4} S1_{16\times 16} =\frac{1}{60} S1_{240\times 240} =115202$

$S1_{240\times 240} :=6912120$

$S2_{240\times 240} :=265427712040$

\bigskip
\noindent
\textbf{$\bullet$ Bimagic square of order 256$\times $256}

\bigskip
$S1_{4\times 4} :=\frac{1}{4} S1_{16\times 16} =\frac{1}{16} S1_{64\times 64} =\frac{1}{64} S1_{256\times 256} =131074$

$S1_{64\times 64} :=2097184$

$S2_{64\times 64} :=91628066144$

$S1_{256\times 560} :=8388736$

$S2_{256\times 256} :=366512264576$

\bigskip
\noindent
\textbf{$\bullet$ Bimagic square of order 343$\times $343}

\bigskip
$S1_{7\times 7} :=\frac{1}{7} S1_{49\times 49} =\frac{1}{49} S1_{343\times 343} =41775$

$S1_{343\times 343} :=20176975$

$S2_{343\times 343} :=1582540680175$

\bigskip
\noindent
\textbf{$\bullet$ Bimagic square of order 400$\times $400}

\bigskip
$S1_{4\times 4} :=\frac{1}{4} S1_{16\times 16} =\frac{1}{20} S1_{80\times 80} =\frac{1}{100} S1_{400\times 400} =320002$

$S1_{400\times 400} :=32000200$

$S2_{400\times 400} :=3413365333400$

\bigskip
\noindent
\textbf{$\bullet$ Bimagic square of order 512$\times $512}

\bigskip
$S1_{8\times 8} :=\frac{1}{4} S1_{32\times 32} =\frac{1}{64} S1_{512\times 512} =1048580$

or

$S1_{4\times 4} :=\frac{1}{4} S1_{16\times 16} =\frac{1}{32} S1_{128\times 128} =\frac{1}{128} S1_{512\times 512} =524290$

$S1_{512\times 512} :=67109120$

$S2_{512\times 512} :=11728191138560$

\bigskip
\noindent
\textbf{$\bullet$ Bimagic square of order 600$\times $600}

\bigskip
$S1_{8\times 8} :=\frac{1}{15} S1_{120\times 120} =\frac{1}{75} S1_{600\times 600} =1440004$

$S1_{600\times 600} :=1018000300$

$S2_{600\times 600} :=25920108000100$

\bigskip
\noindent
\textbf{$\bullet$ Bimagic square of order 625$\times $625}

\bigskip
$S1_{5\times 5} :=\frac{1}{5} S1_{25\times 25} =\frac{1}{25} S1_{125\times 125} =\frac{1}{125} S1_{625\times 625} =976565$

$S1_{125\times 125} :=24414125$

$S2_{125\times 125} :=6357853190125$

$S1_{625\times 625} :=1220706625$

$S2_{625\times 625} :=31789265950625$

\bigskip
\noindent
\textbf{$\bullet$ Bimagic square of order 729$\times $729}

\bigskip
$S1_{9\times 9} :=\frac{1}{81} S1_{729\times 729} =2391489$

$S1_{729\times 729} :=193710609$

$S2_{729\times 729} :=68630571075249$

\bigskip
\noindent
\textbf{$\bullet$ Bimagic square of order 800$\times $800}

\bigskip
$S1_{4\times 4} :=\frac{1}{10} S1_{160\times 160} =\frac{1}{50} S1_{800\times 800} =5120008$

$S1_{800\times 800} :=256000400$

$S2_{800\times 800} :=109226922666800$

\bigskip
\noindent
\textbf{$\bullet$ Bimagic square of order 960$\times $960}

\bigskip
$S1_{4\times 4} :=\frac{1}{4} S1_{16\times 16} =\frac{1}{60} S1_{240\times 240} =\frac{1}{240} S1_{960\times 960} =1843202$

$S1_{960\times 960} :=442368480$

$S2_{960\times 960} :=271791341568160$

\bigskip
\noindent
\textbf{$\bullet$ Bimagic square of order 1000$\times $1000}

\bigskip
$S1_{8\times 8} :=\frac{1}{5} S1_{40\times 40} =\frac{1}{25} S1_{200\times 200} =\frac{1}{125} S1_{1000\times 1000} =4000004$

$S1_{1000\times 1000} :=500000500$

$S2_{1000\times 1000} :=333333833333500$

\bigskip
\noindent
\textbf{$\bullet$ Bimagic square of order 1024$\times $1024}

\bigskip
$S1_{4\times 4} :=\frac{1}{4} S1_{16\times 16} =\frac{1}{16} S1_{64\times 64} =\frac{1}{64} S1_{256\times 256} =\frac{1}{256} S1_{1024\times 1024} =2097154$

$S1_{1024\times 1024} :=536871424$

$S2_{1024\times 1024} :=375300505818624$

\bigskip
\noindent
\textbf{$\bullet$ Bimagic square of order 1200$\times $1200}

\bigskip
$S1_{4\times 4} :=\frac{1}{4} S1_{16\times 16} =\frac{1}{60} S1_{240\times 240} =\frac{1}{300} S1_{1200\times 1200} =2880002$

$S1_{1200\times 1200} :=864000600$

$S2_{1200\times 1200} :=829440864000200$

\bigskip
\noindent
\textbf{$\bullet$ Bimagic square of order 1331$\times $1331}

\bigskip
$S1_{11\times 11} :=\frac{1}{11} S1_{121\times 121} =\frac{1}{121} S1_{1331\times 1331} =9743591$

$S1_{1331\times 1331} :=1178974511$

$S2_{1331\times 1331} :=1392417235445951$

\bigskip
\noindent
\textbf{$\bullet$ Bimagic square of order 1600$\times $1600}

\bigskip
$S1_{4\times 4} :=\frac{1}{4} S1_{16\times 16} =\frac{1}{20} S1_{80\times 80} =\frac{1}{100} S1_{400\times 400} =\frac{1}{400} S1_{1600\times 1600} =5120002$

$S1_{1600\times 1600} :=2048000800$

$S2_{1600\times 1600} :=3495255381333600$

\bigskip
\noindent
\textbf{$\bullet$ Bimagic square of order 2000$\times $2000}

\bigskip
$S1_{4\times 4} :=\frac{1}{4} S1_{16\times 16} =\frac{1}{20} S1_{80\times 80} =\frac{1}{100} S1_{400\times 400} =\frac{1}{500} S1_{2000\times 2000} =8000002$

$S1_{2000\times 2000} :=4000001000$

$S2_{2000\times 2000} :=106667066667000$

\bigskip
\noindent
\textbf{$\bullet$ Bimagic square of order 2048$\times $2048}

\bigskip
$S1_{8\times 8} :=\frac{1}{4} S1_{32\times 32} =\frac{1}{64} S1_{512\times 512} =\frac{1}{256} S1_{2048\times 2048} =16777220$

or

$S1_{4\times 4} :=\frac{1}{4} S1_{16\times 16} =\frac{1}{32} S1_{128\times 128} =\frac{1}{128} S1_{512\times 512} =\frac{1}{512} S1_{2048\times 2048}=8388610$

$S1_{2048\times 2048} :=4294968320$

$S2_{2048\times 2048} :=12009603301288960$

\bigskip
\noindent
\textbf{$\bullet$ Bimagic square of order 2400$\times $2400}

\bigskip
$S1_{8\times 8} :=\frac{1}{15} S1_{120\times 120} =\frac{1}{75} S1_{600\times 600} =\frac{1}{300} S1_{2400\times 2400} =23040004$

$S1_{2400\times 2400} :=6912001200$

$S2_{2400\times 2400} :=26542086912000400$

\bigskip
\noindent
\textbf{$\bullet$ Bimagic square of order 2401$\times $2401}

\bigskip
$S1_{7\times 7} :=\frac{1}{7} S1_{49\times 49} =\frac{1}{49} S1_{343\times 343} =\frac{1}{343} S1_{2401\times 2401} =20176807$

$S1_{343\times 343} :=141237649$

$S2_{343\times 343} :={\rm 379963271712114}0$

$S1_{2401\times 2401} :=6920644801$

$S2_{2401\times 2401} :=26597429019848001$

\bigskip
\noindent
\textbf{$\bullet$ Bimagic square of order 3000$\times $3000}

\bigskip
$S1_{8\times 8} :=\frac{1}{15} S1_{120\times 120} =\frac{1}{75} S1_{600\times 600} =\frac{1}{375} S1_{3000\times 3000} =36000004$

$S1_{3000\times 3000} :=13500001500$

$S2_{3000\times 3000} :=81000013500000500$

\bigskip
\noindent
\textbf{$\bullet$ Bimagic square of order 3200$\times $3200}

\bigskip
$S1_{4\times 4} :=\frac{1}{10} S1_{160\times 160} =\frac{1}{50} S1_{800\times 800} =1280002$

$S1_{3200\times 3200} :=256000400$

$S2_{3200\times 3200} :=109226922666800$

\bigskip
\noindent
\textbf{$\bullet$  Bimagic square of order 4000$\times $4000}

\bigskip
$S1_{4\times 4} :=\frac{1}{10} S1_{160\times 160} =\frac{1}{50} S1_{800\times 800} =\frac{1}{250} S1_{4000\times 4000} =128000008$

$S1_{4000\times 4000} :=32000002000$

$S2_{4000\times 4000} :=341333365333334000$

\bigskip
\noindent
\textbf{$\bullet$ Bimagic square of order 4096$\times $4096}

\bigskip
$S1_{4\times 4} :=\frac{1}{4} S1_{16\times 16} =\frac{1}{16} S1_{64\times 64} =\frac{1}{64} S1_{256\times 256} =\frac{1}{256} S1_{1024\times 1024} =\frac{1}{1024} S1_{4096\times 4096}=33554434$

$S1_{4096\times 4096} :=34359740416$

$S2_{4096\times 4096} :=384307202562021376$

\begin{center}
---------------------------
\end{center}

\end{document}